\documentclass[letterpaper, 10 pt, conference]{ieeeconf}  

\IEEEoverridecommandlockouts                              

\usepackage{amsmath} 
\usepackage{amssymb,amsthm}  

\usepackage{colortbl}

\newtheorem{theorem}{Theorem}[section]
\newtheorem{lemma}{Lemma}[section]
\newtheorem{proposition}{Proposition}[section]
\newtheorem{cor}{Corollary}[section]

\theoremstyle{definition}

\newtheorem{assumption}{Assumption}[section]

\theoremstyle{remark}

\newcommand{\norm}[1]{\|#1\|}

\newcommand{\scal}[2]{\left\langle{#1},{#2}\right\rangle}

\newcommand{\R}{\mathbb{R}}
\newcommand{\N}{\mathbb{N}}
\newcommand{\Db}{\mathcal{D}_{\delta}}
\newcommand{\Dex}{\mathcal{D}}
\newcommand{\xs}{\bar{x}}
\newcommand{\ys}{\bar{y}}
\newcommand{\argmin}{\mathrm{argmin}}

\title{\LARGE \bf
Regularization properties of dual subgradient flow
}

\author{Vassilis Apidopoulos$^{1}$, Cesare Molinari$^{2}$, Lorenzo Rosasco$^{1,3,4}$ and Silvia Villa$^{2}$
\thanks{}
\thanks{$^{1}$MaLGa Center, DIBRIS, Università di Genova, Genoa, Italy
        {\tt\small vassilis.apid@gmail.com}}%
\thanks{$^{2}$MaLGa Center, DIMA, Università di Genova, Genoa, Italy
        {\tt\small cecio.molinari@gmail.com} and {\tt\small villa@dima.unige.it}}%
\thanks{$^{3}$ CBMM, Massachusets Institute of Technology, Cambridge, MA, USA}
\thanks{$^{4}$IIT, Genoa, Italy
        {\tt\small lorenzo.rosasco@unige.it}}%
}

\begin{document}

\maketitle
\thispagestyle{empty}
\pagestyle{empty}

\begin{abstract}

Dual gradient descent combined with early stopping represents an efficient alternative to the Tikhonov variational approach when the regularizer is strongly convex. However, for many relevant applications, it is crucial to deal with regularizers which are only convex. In this setting, the dual problem is non smooth, and dual gradient descent cannot be used. In this paper, we study the regularization properties of a subgradient dual flow, and we show that the proposed procedure achieves the same recovery accuracy as penalization methods, while being more efficient from the computational perspective.

\end{abstract}

\section{Introduction}
We study a continuous dynamical system to recover an
unknown signal $x_*\in\R^p$ from noisy linear measurements 
\begin{equation}
\label{eq:prob}
f^\delta=Ax_*+\epsilon.    
\end{equation}
In this equation, $f^\delta \in\mathbb{R}^d$ is the  vector of measurements, $A\colon\mathbb{R}^p\to\mathbb{R}^d$ is a linear operator and $\epsilon\in\mathbb{R}^d$ is the noise. We assume that the noise is bounded, namely that $\|\epsilon\|\leq \delta$. Moreover, we consider the case $d<p$ and from the modeling point of view we assume that among all the solutions of the equation $Ax_*=Ax$, the vector $x_*$ minimizes a given  function $J$. The main example we have in mind is sparsity, where $J$ is the $\ell_1$ norm. Such problems have been widely studied in machine learning \cite{hastie2015statistical}, inverse problems \cite{engl1996regularization,benning2018modern}, and signal processing/compressed sensing \cite{candes2006robust,foucart2013invitation}. The most classical approach is based on explicit (Tikhonov/variational) regularization, which is obtained by the minimization of the weighed sum of two terms: a penalization of the equality constraints $Ax=f^\delta$ (data fit term); and the regularizer $J$. This requires the selection of an appropriate trade-off parameter, which is typically the computational bottleneck of the regularization approach. Motivated by the stability and the reconstruction properties of first order methods in machine learning, in this paper we follow a different route, based on the so-called {\it implicit regularization} (iterative regularization in the inverse problems literature). Such an approach combines regularization properties with computational efficiency, since it is based on the early stopping of an iterative procedure used to minimize $J$ among the solutions of \eqref{eq:prob}. The most well-known example is gradient descent, a.k.a. Landweber iteration \cite{engl1996regularization}. Several extensions of this basic approach have been proposed to deal with more general regularizations, which are strongly convex \cite{matet2017don}, or just convex \cite{molinari2021iterative,molinari2022iterative}. These regularization methods depend on the choice of a suitable discrete algorithm or continuous dynamical system. In this paper we focus on the convex case, and we propose to find a stable reconstruction of $x_*$ by early stopping a subgradient flow in the dual space $\mathbb{R}^d$. We prove that such a procedure has the same reconstruction properties of Tikhonov regularized solutions in terms of dependence from the error, while being cheaper from the computational perspective. This problem has been studied already in \cite{burger2004convergence,burger2007error,boct2021convergence}. We leave the comments and the comparison of our work with the previous literature after the main results.\\
The paper is organized as follows. In Section \ref{sectionprelim} we present the main problem, the setting considered and the precise assumptions on the objects involved. In Section \ref{secdualflow} we present the dual subgradient flow and the main theoretical results about iterative regularization through early-stopping. Some comments on the assumptions, on the dynamical system and on the results are reported in Section \ref{sectiondiscussion}. In Section \ref{sectionconvergence} we provide the convergence and stability analysis for the dual subgradient flow. We conclude with a summary of the work done and some possible future developments in Section \ref{sectionconcl}.

\section{Problem and assumptions}\label{sectionprelim}
We assume that the unknown signal $x_*$ is a solution of the constrained optimization problem
\begin{equation}\label{minproblem}
    \underset{x\in \R^p}{\min} \quad J(x) ~ : ~~ Ax=f,
\end{equation}
where $f:=Ax_*$ is the exact datum, typically not available in practice. 
We make the following general assumptions on the objects involved.
\begin{assumption}\label{genass}
The operator $A:\R^p\to \R^d$ is linear,  $J:\R^p\to \bar{\mathbb{R}}:=\mathbb{R}\cup\{+\infty\}$ is a proper, convex and lower semi-continuous function and $f$ is a given vector in $\R^d$. 
\end{assumption}
We illustrate examples of regularizers $J$ which are relevant in applications and satisfy the previous assumptions.
\begin{itemize}
    \item When $J(\cdot)=\norm{\cdot}_{2}^{2}$, problem \eqref{minproblem} is used to find the minimal norm solution to the linear equation $Ax=f$. The unique minimizer is given by the Moore-Penrose pseudo-inverse of $A$ applied to the exact data, namely $A^{\dagger} f$ (see e.g. \cite{engl1996regularization}).
    \item The choice of $J(\cdot)=\norm{\cdot}_{1}$ is well-known to induce sparsity in the solution. Indeed, it is the convex relaxation of the $0-$seminorm, that counts the non-zero elements in a vector \cite{candes2006robust}.
    \item $J(\cdot)=\norm{\cdot}_{TV}$ is the total variation regularizer. It is used in image reconstruction to maintain sharp edges between different zones with constant color \cite{rudin1992nonlinear}.
    \item $J(\cdot)=\norm{\cdot}{}_*$ is the nuclear norm, namely the $1-$norm of the singular values of a matrix. It is an approximation of the rank and so its minimization is used to maximize the linear correlation between rows or columns \cite{cai2010singular}. 
\end{itemize}
Except for the $2-$norm, the examples above involve regularizers  that are neither differentiable nor strongly-convex, but only lower-semicontinuous and convex.

We also introduce the dual problem associated to the exact problem \eqref{minproblem}, namely
\begin{equation}\label{dualprobl}
    \underset{y\in \R^d}{\min} \quad \Dex(y):= J^{\star}(-A^*y) + \langle f, y\rangle,
\end{equation}
where $J^{\star}$ is the Fenchel conjugate of $J$, defined as $J^{\star}(x):=\sup_{z\in \R^p} \ \langle x,z\rangle - J(z)$. We say that a pair $\left(\bar{x},\bar{y}\right)\in\R^p\times \R^d$ is a {\it primal-dual solution} of \eqref{minproblem}-\eqref{dualprobl} if it satisfies the following KKT conditions:
\begin{equation}\label{solutionKKT}
		\begin{cases}
-A^*\bar{y}\in \partial J(\bar{x}) \\
	~	A\bar{x} =f,    
		\end{cases}
\end{equation}
where $\partial J(x)$ is the subdifferential of $J$ at $x\in\R^p$, namely the set-valued operator defined by
$$\partial J(x) := \left\{x'\in\R^p: \ J(z)\geq J(x) + \langle x',z-x\rangle \ \forall z\in\R^p \right\}.$$ 
In the next, we assume the existence of a primal-dual solution.
\begin{assumption}\label{assumptionsolutions}
    There exists an optimal primal-dual pair $(\bar{x},\bar{y})\in\R^p\times\R^d$ for the exact problem, namely satisfying \eqref{solutionKKT}.
\end{assumption}

The previous assumption is classical. Some comments are given in the discussion paragraph \ref{discass}.

Notice that while we are searching to approximate a solution of problem \eqref{minproblem}, the related datum $f$ is typically not available in practice. In this setting, the aim of regularization is to recover, in a stable way with respect to the noise level $\delta$, an approximate solution to \eqref{minproblem} where only a noisy version of $f$, denoted by $f^{\delta}$, is available. As already anticipated in the Introduction, we focus on the worst-case deterministic scenario, in which $\delta \geq 0$ denotes the noise level. More precisely, we make the following assumption on the inexact data.
\begin{assumption}\label{assumptiondelta}
For some $\delta\geq 0$, the element $f^{\delta}\in\R^d$ satisfies 
\begin{equation}\label{bounddelta}
    \norm{f^{\delta}-f} \leq \delta.
\end{equation}      
\end{assumption}
Any regularization method should be such that  when $\delta$ approaches 0, the regularized solution approaches the ideal vector $x_*$. 
In the next section, we will show that the trajectories of a dual dynamical system can be used to build a regularization method.

\section{Dual flow and main results}\label{secdualflow}
In this section we describe the proposed implicit/iterative regularization method to approximate a solution of \eqref{minproblem} in a stable way with respect to the noise level in the data $f$. Then we present the main theoretical results of the paper.

Given $f^{\delta}\in\R^d$ satisfying Assumption~\ref{assumptiondelta}, we define the inexact dual function as
$$\Db(y):= J^{\star}(-A^*y) + \langle f^{\delta}, y\rangle.$$
Consider $y_{0}\in\R^d$ such that $\partial \Db(y_0)\neq \emptyset$. For $t\geq 0$ we are interested in the behavior of the trajectory $\{y(t)\}_{t\geq0}$  generated by the subgradient flow on $\Db$, i.e.: $y(0)=y_{0}$ and, for every $t>0$,
\begin{equation}\label{dualflow}
    \dot{y}(t)+\partial \Db(y(t))\ni 0.
\end{equation}

In the next, we assume a technical condition on the composition of $J^\star$ with the linear operator $A$, that will allow us to re-write the dual subgradient flow as a coupled dynamical system.
\begin{assumption}\label{assumptionchainrule}
    For every $y\in\R^d$, the following chain rule equality holds:
    \begin{equation}\label{chainrule}
        \partial \left[ J^{\ast}(-A^{\ast}\cdot)\right](y)= -A\partial J^{\ast}(-A^{\ast}y).
    \end{equation}
\end{assumption}
The previous assumption is classical. Some comments are given in the discussion paragraph \ref{discass}. By using the chain rule in \eqref{chainrule}, the differential inclusion in \eqref{dualflow} is equivalent to
\begin{equation*}
    \dot{y}(t)\in A\partial J^{\ast}(-A^{\ast}y(t)) -f^{\delta}
\end{equation*}
and so to the following dynamical system:
\begin{equation}\label{DGFP}
    \begin{cases}
        x(t)&\in \partial J^{\ast}(-A^{\ast}y(t)) \\
        \dot{y}(t)&=Ax(t)-f^{\delta}.
    \end{cases}
\end{equation}

More discussion on the dual flow can be found in paragraph \ref{discussflow}.



In order to introduce the main results, we recall the definition of the {\it Bregman divergence} related to $J$. For $(x_{1},x_{2})\in \R^{p}\times\R^p$ and $u\in \partial J(x_{2})$, the Bregman divergence of $J$ between $x_{1}$ and $x_{2}$ is given by
\begin{equation}
    D^{u}(x_{1},x_{2}):=J(x_{1})- J(x_{2}) -\scal{u}{x_{1}-x_{2}}.
\end{equation}
To show that $x(t)$, the trajectory of the dual subgradient flow in \eqref{DGFP}, is indeed approaching a solution of the exact problem \eqref{minproblem}, one option is to show bounds for the classical duality gap (see the definition in \eqref{dualitygap} in the discussion section). In this work, we follow a different strategy, proving rates for both the {\it Lagrangian gap} $D^{-A^{\ast}\bar{y}}(x(t),\bar{x})$ and the {\it feasibility gap} $\norm{Ax(t)-f}$. Indeed, if $D^{-A^{\ast}\bar{y}}(x',\bar{x})=0$ and $Ax'=f$ for some $x'\in\R^p$, then $x'$ is a solution for problem \eqref{minproblem} (see e.g. \cite[Lemma 5]{molinari2021iterative}). So the combination of the two quantities is indeed a measure of (primal) optimality for the point $x'$. We are now ready to give the main results of the paper, concerning the regularization properties of the dual subgradient flow \eqref{DGFP}.



\begin{theorem} [Rates]\label{basiccor}
Let Assumptions \ref{genass}, \ref{assumptionsolutions}, \ref{assumptiondelta} and \ref{assumptionchainrule} hold true and $(x(t),y(t))_{t\geq 0}$ be a solution of \eqref{DGFP}. Consider the tail-average sequence \begin{equation}
    \hat{x}(t):=\frac{2}t\int_{\frac{t}{2}}^{t}x(s)ds.
\end{equation}
Then the following rates of convergence hold true for all $t\geq 0$: for $C_{0}:=\norm{y(0)-\bar{y}}$,
\begin{equation}\label{rateBregman}
    D^{-A^{\ast}\bar{y}}(\hat{x}(t),\bar{x}) \leq \frac{C_{0}^{2}}{t}+2C_{0}\delta +\delta^{2}t;
\end{equation}
\begin{equation}\label{ratefeasibility}
    \norm{A\hat{x}(t)-f}^{2} \leq \frac{2C_{0}^{2}}{t^{2}}+6C_{0}\frac{\delta}{t} +6\delta^{2};
\end{equation}
\begin{equation}\label{estimateDinverse}
          D^{-A^{\ast}y(t)}(\bar{x},x(t)) \leq \frac{C_{0}^{2}}{2t}+\frac{3C_{0}}{2}\delta +\frac{3}{2}\delta^{2}t.
    \end{equation}    
\end{theorem} 

The proof of the Theorem \ref{basiccor} can be found in the related paragraph \ref{subsectionProofteo}. From Theorem \ref{basiccor}, we deduce the following early-stopping result.
\begin{cor}[Early stopping]\label{earlystop}
  Let Assumptions \ref{genass}, \ref{assumptionsolutions}, \ref{assumptiondelta} and \ref{assumptionchainrule} hold true, $(x(t),y(t))_{t\geq 0}$ be a solution of \eqref{DGFP} and set $\hat{x}(t):=\frac{2}t\int_{\frac{t}{2}}^{t}x(s)ds$. Then, for $t^{\ast}:=\frac{1}{\delta}\norm{y(0)-\bar{y}}$, it holds:
  \begin{equation}
      D^{-A^{\ast}\bar{y}}(\hat{x}(t^*),\bar{x}) \leq 4\norm{y(0)-\bar{y}}\delta,
  \end{equation}
  \begin{equation}
      \norm{A\hat{x}(t^*)-f} \leq \sqrt{14}\delta.
  \end{equation}
\end{cor}
For the comments on the previous results, see the discussion paragraph \ref{discussres}.

\section{Discussion}\label{sectiondiscussion}

In this section, we discuss the assumptions used in our analysis; the proposed iterative regularization method, the dual subgradient flow; and the theoretical results obtained. We also compare our results first with the ones for the same dual flow but under stronger assumptions, namely when $J$ is strongly convex; and then with the results given by explicit/variational regularization by Tikhonov path. We conclude with some remarks about possible discretization of the dual subgradient flow, in order to get optimization algorithms.
\subsection{Discussion on the assumptions}\label{discass}
Assumption \ref{assumptionsolutions} and \ref{assumptionchainrule} are classical in saddle-point problems and they hold, for instance, under the so called \textit{qualification conditions} 
(see e.g. \cite[Proposition $3.28$, Theorem 3.30, Corollary 3.31 and Theorem 3.51]{peypouquet2015convex}). In particular, we recall that a sufficient condition for \eqref{chainrule} is the existence of a point $y_0\in\R^d$ such that $J^{\star}$ is continuous at $-A^*y_0$. \\
Regarding Assumption \ref{assumptionsolutions}, note that, if $(\bar{x},\bar{y})$ is a primal-dual solution, in particular $\bar{x}$ is a minimizer for \eqref{minproblem} and $\bar{y}$ for \eqref{dualprobl}. The opposite is not true, in general. We briefly discuss here sufficient conditions in order to ensure existence of a primal-dual solution.
If there exists a (dual) solution $\bar{y}\in\R^d$ to problem \eqref{dualprobl} and the condition in \eqref{chainrule} holds, then Assumption \ref{assumptionsolutions} holds true. Indeed, by generalized Fermat's rule (see, for instance, \cite[Theorem 3.24]{peypouquet2015convex}), $\bar{y}$ satisfies
\begin{align*}
0 & \in \partial \left[J^{\star}(-A^*\cdot)+\langle f, \cdot \rangle \right](\bar{y})\\
& = -A\partial J^{\star}(-A^*\bar{y})+f.
\end{align*}
Denoting by $\bar{x}$ an element in $\partial J^{\star}(-A^*\bar{y})$ for which the previous inclusion is satisfied and using the fact that $(\partial J^{\star})^{-1}=\partial J$ (see, for instance, \cite[Proposition 3.59]{peypouquet2015convex}, we get that the pair $(\bar{x},\bar{y})$ is a primal-dual solution accordingly to the KKT conditions in \eqref{solutionKKT}. \\
On the other hand, if there exists a primal solution $\bar{x}$ and a point $x_0\in\R^p$ such that $Ax_0=f$ (feasibility) and where $J$ is continuous, then there exists a dual solution $\bar{y}\in\R^d$ such that $(\bar{x},\bar{y})$ is an optimal primal-dual pair. Indeed, it is easy to show that, under the general Assumption \ref{genass}, the following chain rule holds:
$$\partial \left[\iota_{\left\{f\right\}}\circ A\right](\cdot)=A^*\partial \iota_{\left\{f\right\}} (A\cdot),$$
where $\iota_{\left\{f\right\}}$ denotes the classical indicator function of convex analysis; namely, $\iota_{\left\{f\right\}}(y)$ is $0$ if $y=f$ and $+\infty$ otherwise. Moreover, at the point $x_0$, the function $\iota_{\left\{f\right\}}\circ A$ assumes a bounded value and $J$ is continuous. So also the sum rule for subdifferentials holds (see, for instance, \cite[Theorem 3.30]{peypouquet2015convex}). Then the primal solution $\bar{x}$ satisfies $A\bar{x}=f$ and
\begin{align*}
0& \in \partial \left[ J+\iota_{\left\{f\right\}}\circ A\right](\bar{x})\\
& =     \partial J(\bar{x})+A^*\partial\iota_{\left\{f\right\}}(A\bar{x})\\
& =     \partial J(\bar{x})+A^*\R^d.
\end{align*}
The previous means precisely that there exists a point $\bar{y}\in\R^d$ such that $(\bar{x},\bar{y})$ satisfies the KKT conditions in \eqref{solutionKKT}.
\\
\subsection{Discussion on the dual flow}\label{discussflow}
If $f^{\delta}$ belongs to the range of operator $A$, a natural approach to deal with problem \eqref{minproblem} is to consider the subgradient flow on the (noisy) primal objective function, i.e. on $x \mapsto J(x) + \iota_{\{f^\delta\}}(Ax)$. Under standard qualification conditions, 
for $t\geq 0$, the primal flow reads as
\begin{equation}\label{primal flow}
    \dot{x}(t) +\partial J(x(t)) +A^{\ast}\partial\iota_{\{f^\delta\}}(Ax(t))\ni 0.
\end{equation}
Nevertheless, we shall stress out that this system forces the trajectory to be feasible accordingly to the inexact constraint (namely, $Ax(t)=f^{\delta}$ for every $t \geq 0$), which is not desirable in general. Furthermore, the assumption $f^{\delta}\in \text{Im}(A)$ is not necessary in our setting for considering the dual flow \eqref{dualflow}, contrarily to the primal one \eqref{primal flow}. These are the main motivations for studying the dual flow \eqref{dualflow} instead of \eqref{primal flow}. 
\ \\
Notice that, if the function $J$ is twice differentiable with invertible Hessian, the system \eqref{DGFP} can be also written as 
\begin{equation}\label{DGFP2}
\dot{x}(t) = \left[\nabla^2 J (x(t))\right]^{-1} A^*\left(f^{\delta}-Ax(t)\right).
\end{equation}
So, under strong regularity assumptions that we do not suppose in this work, the system \eqref{DGFP} may be seen as a generalized mirror flow associated to the least-squares problem $Ax=f^{\delta}$, with mirror function given by $J$.

\subsection{Discussion on the results}\label{discussres}

In order to explain the results in Theorem \ref{basiccor} and Corollary \ref{earlystop}, we first recall the definition of the Lagrangian functional related with \eqref{minproblem}:
\begin{equation*}
    \mathcal{L}(x,y) := J(x)+\scal{y}{A{x}-f}.
\end{equation*}
Note that $(\bar{x},\bar{y})$ is a primal-dual solution accordingly to \eqref{solutionKKT} if and only if it is a saddle-point for the Lagrangian; namely if and only if, for every $(x,y)\in\R^p\times\R^d$,
\begin{equation}\label{saddlepoint}
        \mathcal{L}(\bar{x},y)\leq \mathcal{L}(\bar{x},\bar{y}) \leq \mathcal{L}(x,\bar{y}).
\end{equation}
Then the quantity 
\begin{equation}\label{dualitygap}
\sup_{x\in\R^p,\ y\in\R^d} \ \mathcal{L}(x',y)-\mathcal{L}(x,y'),  
\end{equation}
called {\it duality gap}, is a reasonable measure of optimality for the point $(x',y')\in\R^p\times\R^d$; see, for instance, \cite{chambolle2011first}. From \eqref{saddlepoint}, another interesting quantity to measure the optimality of $x'\in\R^p$ is 
\begin{equation}\label{lagrgap}
\mathcal{L}(x',\bar{y})-\mathcal{L}(\bar{x},\bar{y}),   
\end{equation}
where $(\bar{x},\bar{y})$ is a given primal-dual solution. The term {\it Lagrangian gap}, used above for $D^{-A^*\bar{y}} (x',\bar{x})$, comes from the fact that it corresponds to \eqref{lagrgap}: indeed, as $A\bar{x}=f$, we have
\begin{equation*}
	\begin{aligned}
		D^{-A^*\bar{y}} (x',\bar{x}) & = J(x')-J(\bar{x})-\langle -A^*\bar{y}, x'-\bar{x} \rangle\\
        & = \mathcal{L}(x',\bar{y})-\mathcal{L}(\bar{x},\bar{y}).
	\end{aligned}
\end{equation*}



In \cite{burger2007error}, the authors proved similar results to the ones in Corollary \ref{earlystop}, but only for $D^{-A^*y(t)}(\bar{x},x(t))$. The previous quantity is not directly related to optimality of the primal variable $x(t)$ if the function $J$ is not strongly-convex. The main difference between our result and the one in \cite{burger2007error} is that we are able to control two different quantities: the {\it Lagrangian gap} $D^{-A^*\ys}(\cdot,\xs)$, with respect to a primal-dual solution pair $(\xs,\ys)$; and the {\it feasibility gap} $\norm{A\cdot-f}$, with respect to the exact data $f$. As already remarked in the previous section, the combination of these two quantities is a good measure of optimality. Indeed, when they are both zero at $x'\in\R^p$, then $x'$ is a solution of the primal optimization problem \eqref{minproblem}; see \cite{molinari2021iterative,molinari2022iterative}. The bounds that we derive show that the dual subgradient flow \eqref{DGFP} on the inexact data $f^\delta$, if opportunely early-stopped, can be used to implicitly regularize the solution of a linear inverse problem affected by worst-case deterministic noise. More precisely, the upper-bounds in Theorem \ref{basiccor} on the Lagrangian gap and the feasibility are composed by the sum of two different terms. The first term is related with the {\it optimization} procedure and so it is vanishing for $t\to +\infty$. The second one is related with {\it stability} with respect to the error on the data. This term is bounded of order $\delta$ for the feasibility and it is diverging for $t\to +\infty$ as $\delta t$ in the case of the Lagrangian gap. Minimizing the upper-bound with respect to $t$ suggests the choice of an optimal stopping time $t^* = C/\delta$ for some constant $C>0$, as shown in Corollary \ref{earlystop}. For this choice, both the quantities $D^{-A^*\ys}(\hat{x}(t^*),\xs)$ and $\norm{A\hat{x}(t^*)-f}$ are of order $\delta$, where $\hat{x}$ is the tail averaged trajectory. This reflects the idea of stability by iterative regularization: an error of order $\delta$ on the data is reflected to an error of order $\delta$ on the solution, without amplifications. 
\subsection{Comparison with results for strongly convex}
When $J$ is strongly convex with constant $\mu>0$, $J^{\star}$ is smooth (see, for instance, \cite[Proposition 14.2]{bauschke2011convex}); meaning that $J^{\star}$ is differentiable with Lipschitz continuous gradient. Then, the differential inclusion in \eqref{DGFP} becomes the classical differential equation
\begin{equation*}
    \begin{cases}
        x(t)&= \ \nabla J^{\ast}(-A^{\ast}y(t)) \\
        \dot{y}(t)&= \ Ax(t)-f^{\delta}.
    \end{cases}
\end{equation*}
Moreover, the Bregman divergence verifies 
$$D^{-A^*\bar{y}}(x,\bar{x}) \geq \frac{\mu}{2}\norm{x-\bar{x}}^2$$
and so the Lagrangian gap controls the distance to the solution.
In particular, our bounds on the Lagrangian gap in Theorem \ref{basiccor} and Corollaries \ref{earlystop} automatically transfer to the quantity $\norm{\hat{x}(t)-\bar{x}}$. It can be shown that similar bounds holds indeed for the non-averaged iterate $x(t)$. See \cite{engl1996regularization,burger2004convergence}, the recent work \cite{DGF22} and, for the discrete case, \cite{matet2017don}.
\subsection{Comparison with Tikhonov path}

Under Assumptions \ref{genass}, \ref{assumptionsolutions} and \ref{assumptiondelta}, the Tikhonov regularization path for $s\geq 0$ is defined by
\begin{equation}
    u(s) \in \argmin_{u\in\R^p} \left\{ J(u)+\frac{s}{2}\norm{Au-f^{\delta}}^2\right\}
\end{equation}
or, equivalently, by
\begin{equation}\label{ee}
\begin{split}
0&\in \partial \left( J(\cdot) + \frac{s}{2}\norm{A\cdot-f}^2\right)\left(u(s)\right)\\
& = \partial J(u(s)) + s A^*\left(Au(s)-f^\delta\right).
\end{split}
\end{equation}
Denoting $v(s):=Au(s)-f^\delta$, we have that $-sA^* v(s)\in \partial J(u(s))$ and so that $(u(s),v(s))$ satisfies the system
\begin{equation}\label{e}
\begin{cases}
u(s) \ \in & \partial J^\star (-sA^*v(s));\\
v(s)\ =& Au(s)-f^\delta.
\end{cases}
\end{equation}
Assuming enough regularity on the functional $J$ and differentiating \eqref{ee} with respect to the variable $s$, we get $0=\nabla^2 J(u(s)) \dot{u}(s) + A^*\left(Au(s)-f^{\delta}\right)+sA^*A\dot{u}(s)$ and
\begin{equation}\label{f2}
\dot{u}(s) = \left[\nabla^2 J (u(s)) + s A^*A\right]^{-1} A^*\left(f^{\delta}-Au(s)\right).
\end{equation}
Notice the similarity between \eqref{DGFP}-\eqref{DGFP2} and \eqref{e}-\eqref{f2}. While \eqref{DGFP} is a dynamical inclusion, the inclusion in \eqref{e} does not involve any derivative.

Analogously to Theorem \ref{basiccor} and Corollary \ref{earlystop}, for the Tikhonov path we have the following bounds and regularization properties.
\begin{theorem} \label{teotikh}
Under Assumptions \ref{genass}, \ref{assumptionsolutions} and \ref{assumptiondelta}, let $s>0$ and $u(s)$ be the Tikhonov path defined in \eqref{ee} and $v(s)=Au(s)-f^{\delta}$. Then the following bounds hold true:
\begin{equation}\label{DsymTikh}
    D^{-A^{\ast}\bar{y}}(u(s),\bar{x}) +D^{-sA^{\ast}v(s)}(\bar{x},u(s)) \leq \frac{1}{4s}\left(\norm{\bar{y}}+s\delta\right)^{2};
\end{equation}	
\begin{equation}\label{feasibilityTykh}
    \norm{Au(s)-f} \leq \frac{\norm{\bar{y}}}{s}+\delta.
\end{equation}	
In particular, for $s^*:=\norm{\bar{y}} / \delta$,
	$$D^{-A^{\ast}\bar{y}}(u(s^*),\bar{x}) +D^{-sA^{\ast}v(s^*)}(\bar{x},u(s^*)) \leq \norm{\bar{y}} \delta$$
	and
	$$\norm{Au(s^*)-f} \leq 2\delta.$$

\end{theorem} 
The proof of Theorem \ref{teotikh} can be  seen as a straightforward generalization of classic results (see for instance \cite{engl1996regularization,burger2004convergence})  and is given for completeness in Section \ref{TYKgen}.
Note the similarities and differences between the results for the iterative (implicit) regularization properties given by the dual subgradient flow in  Theorem \ref{basiccor} and Corollary \ref{earlystop} and the ones for the Tikhonov (explicit) regularization in Theorem \ref{teotikh}. The main difference is that for Tikhonov path the results hold for the trajectory $u(s)$ without averaging. 
On the other hand, the final results are qualitatively very similar. Indeed, in the two cases the best early-stopping $t^*$ or $s^*$ is proportional to $1/\delta$ and, at that value, both the Lagrangian gap and the feasibility achieve a bound of order $\delta$.

\subsection{Discretization}\label{secdiscr}
There are two possibilities in order to obtain algorithms by time discretization of the flow defined in \eqref{DGFP}. The first option is a backward-implicit discretization, that leads to the following scheme: given $y_0 \in \R^d$, for every $k\in\N$,
\begin{equation*}
    \min_{y\in\R^d} \left\{ J^{\ast}(-A^*y)+\langle f^\delta, y\rangle+ \frac{1}{2\gamma}\norm{y-y_k}^2\right\};
\end{equation*}
namely, the proximal-point algoritm on the inexact dual function $\Db$. Another way to see the implicit discretization, on the primal variable, is the following: consider $-A^{\ast}y_{k+1}\in \partial J(x_{k+1}),$
where 
\begin{equation*}
        -A^*y_{k+1}=-A^*y_k-\gamma A^*\left(Ax_{k+1}-f^{\delta}\right).
\end{equation*}
The relations above imply that
$0 \in \partial J(x_{k+1}) + \gamma A^*A x_{k+1} -\gamma A^*f^{\delta}+A^*y_k$.
The previous means that $x_{k+1}$ is a solution to the optimization problem 
\begin{equation*}
    \min_{x\in\R^p} \left\{ J(x)+ \frac{\gamma}{2}\norm{Ax-\left(f^{\delta}-\frac{y_k}{\gamma}\right)}^2\right\}.
\end{equation*}
Note that also in this case each iterative step is implicit and so computationally expensive. This is also known as Bregman iteration and has been studied in \cite{burger2007error}.

On the other hand, a forward-explicit discretization of \eqref{DGFP} leads to a subgradient-type algorithm of the following form: pick $x_k\in \partial J^{\ast}(-A^{\ast}y_k)$ and update the dual variable by
\begin{equation}\label{subgradalg}
       y_{k+1}=y_k+\gamma \left(Ax_k-f^{\delta}\right).
\end{equation} 
We leave the study of the regularization properties of this algorithm for future work. We remark here that in \cite{molinari2021iterative} and \cite{molinari2022iterative} another method, namely the primal-dual algorithm \cite{chambolle2011first,condat2013primal,vu2013splitting,rasch2020inexact}, has been studied in order to perform iterative regularization for problem \eqref{minproblem}.

\section{Convergence analysis}\label{sectionconvergence}

In this section we provide the proof of the main Theorem \ref{basiccor}, which relies in Lyapunov-based techniques. In particular we first recall the notion of a solution to \eqref{DGFP} and then provide a suitable Lyapunov candidate for the system \eqref{DGFP}, which allows to deduce the estimates on the Bregman divergence and the feasibility gap as stated in Theorem \ref{basiccor}.

\subsection{Existence and regularity of solutions}
Let us first recall some basic facts and  properties concerning the existence and the regularity of a solution of \eqref{dualflow} (equivalently \eqref{DGFP}).

\begin{proposition}\label{propexistence}\cite[Theorems $3.1$ \& $3.2$]{brezis1973ope}
    Consider the system \eqref{DGFP}. Then there exists a unique function $y:[0,+\infty)\to\R^d$  such that:
\begin{enumerate}
    \item $y(t)\in \text{dom}\partial \Db$ for all $t\geq 0$.
    \item $y$ is Lipschitz continuous in $[0,+\infty)$.
    \item $x(t) \in \partial J^{\ast}(-A^{\ast}y(t))$ and $\dot{y}(t) =Ax(t)-f^{\delta}$ a.e. in $[0,+\infty)$.
    \item The function $\mathcal{D}_{\delta}\circ y:[0,+\infty)\to\bar{\mathbb{R}}$ is convex, non-increasing and absolutely continuous in $[0+\infty)$, with
    \begin{equation}
        \frac{d}{dt}\mathcal{D}_{\delta}(y(t))= -\norm{Ax(t)-f^{\delta}}^{2} ~ \text{ a.e. in }[0,+\infty)
    \end{equation}
\end{enumerate}
\end{proposition}

\subsection{Lyapunov estimates}
In this paragraph we present the main Lyapunov analysis that will be used for the proof of Theorem \ref{basiccor}.

We first give some preliminary estimates on the dual distance $\norm{y(t)-\bar{y}}$.
\begin{lemma}\label{lemma bound y}
 Let $(x(t),y(t))_{t\geq 0}$ be a (coupled) solution of \eqref{DGFP}. Then for all $t\geq0$, it holds:
     \begin{equation}\label{boundy}
         \norm{y(t)-\bar{y}} \leq \norm{y(0)-\bar{y}} +\delta t
     \end{equation} 
     \begin{equation}\label{boundDsym}
     \begin{aligned}
         \int_{0}^{t}&\left(D^{-A^{\ast}\bar{y}}(x(s),\bar{x})  +D^{-A^{\ast}y(s)}(\bar{x},x(s)) \right) ds \leq \\
         & \quad \qquad \leq \frac{\norm{y(0)-\bar{y}}^{2}}{2} + \norm{y(0)-\bar{y}}\delta t + \frac{(\delta t)^2}{2}       
     \end{aligned}
     \end{equation} 
\end{lemma}


{\em Proof.}
By taking the derivative of $\norm{y(t)-\bar{y}}^{2}$ with respect to time and using \eqref{DGFP}, for all $t\geq 0$, we have: 
\begin{equation}\label{derivative y}
\begin{aligned}
\frac{1}{2}\frac{d}{dt} \norm{&y(t)-\bar{y}}^2 =\langle \dot{y}(t), y(t)-\bar{y} \rangle\\
&=\langle Ax(t)-f^{\delta}, y(t)-\bar{y}\rangle\\
&=  \langle Ax(t)-f, y(t)-\bar{y}\rangle + \langle f-f^{\delta}, y(t)-\bar{y}\rangle \\
&=\langle x(t)-\bar{x}, A^*y(t)-A^*\bar{y}\rangle +\langle f-f^{\delta}, y(t)-\bar{y}\rangle \\
&=-\langle x(t)-\bar{x}, \left(-A^*y(t)\right)-\left(-A^*\bar{y}\right)\rangle \\ & \quad +\langle f-f^{\delta}, y(t)-\bar{y}\rangle \\
&=-\left(D^{-A^{\ast}\bar{y}}(x(t),\bar{x})  +D^{-A^{\ast}y(t)}(\bar{x},x(t)) \right) \\ 
& \quad+\langle f-f^{\delta}, y(t)-\bar{y}\rangle
\end{aligned}
\end{equation}
By applying the Cauchy-Schwarz inequality in \eqref{derivative y}, we find: 
\begin{equation}\label{eq: y before bihari2}
\begin{aligned}
  \frac{1}{2}\frac{d}{dt} \norm{y(t)-\bar{y}}^2  &\leq -D^{-A^{\ast}\bar{y}}(x(t),\bar{x})  \\ & \quad -D^{-A^{\ast}y(t)}(\bar{x},x(t))  +\delta\norm{y(t)-\bar{y}}    
\end{aligned}
\end{equation}
thus by integrating \eqref{eq: y before bihari2} on $[0,t]$, we have:
\begin{equation}\label{eq: y before bihari}
    \frac{1}{2}\norm{y(t)-\bar{y}}^2 \leq \frac{1}{2}\norm{y(0)-\bar{y}}^2 + \int_{0}^{t}\delta \norm{y(s)-\bar{y}} ds
\end{equation}
thus, by applying Lemma \ref{lemma Bihari} for \eqref{eq: y before bihari}, for all $t\geq0$, we obtain:
\begin{equation}\label{boundy2}
    \norm{y(t)-\bar{y}} \leq \norm{y(0)-\bar{y}} +\delta t
\end{equation}
which proves \eqref{boundy}.

By using the estimate \eqref{boundy2} in \eqref{eq: y before bihari}, we find:
\begin{equation}
\begin{aligned}
  & \frac{\norm{y(t)-\bar{y}}^{2}}{2} + \\ & +  \int_{0}^{t}\left(D^{-A^{\ast}\bar{y}}(x(s),\bar{x})  +D^{-A^{\ast}y(s)}(\bar{x},x(s)) \right) ds  \leq \\ & \qquad\qquad \qquad \leq \frac{\norm{y(0)-\bar{y}}^{2}}{2}   +\left(\norm{y(0)-\bar{y}} +\frac{\delta t}{2}\right)\delta t 
\end{aligned}
\end{equation}
 which allows to conclude with the estimate \eqref{boundDsym}
\qed

Here we shall notice that the estimate \eqref{rateBregman} for the averaged tail trajectory $\hat{x}(t)=\int_{\frac{t}{2}}^{t}x(s)ds$ in Theorem \ref{basiccor}  can be directly provided from \eqref{boundDsym} in Lemma \ref{lemma bound y} and the dissipativity property of the energy $\frac{1}{2}\norm{y(t)-\bar{y}}^{2}$ . However, as discussed previously, in order to ensure the regularising property of the primal trajectory $x(t)$ of the system \eqref{DGFP} towards a solution $\bar{x}$, we also need to control the feasibility gap (i.e. $\norm{Ax(t)-f}$). For this purpose we will use a different Lyapunov energy function.
In particular, for all $t\geq 0$, we consider the following function:
\begin{equation}\label{energy definition}
\begin{aligned}
    V(t)&=\frac{1}{2}\norm{y(t)-\bar{y}}^{2} +\int_{0}^{t}s\norm{Ax(s)-f}^{2} ds   \\
    &\quad -\delta\int_{0}^{t}\norm{y(s)-\bar{y}} ds -\delta\int_{0}^{t}s\norm{Ax(s)-f} ds \\ &\quad +\int_{0}^{t} D^{-A^{\ast}\bar{y}}(x(s),\bar{x}) ds+tD^{-A^{\ast}y(t)}\left(\bar{x},x(t)\right)
\end{aligned}
\end{equation}

The next proposition shows that the function $V$ is a Lyapunov function for the system \eqref{DGFP}.

\begin{proposition}\label{propVdecreasing}
    Let $(x(t),y(t))_{t\geq 0}$ be a (coupled) solution of \eqref{DGFP} and $V(t)$ as defined in \eqref{energy definition}. Then the function $V$ is absolutely continuous and non-increasing in $[0,+\infty)$. 
\end{proposition}

{\em Proof.}
Let us first show that $V$ is absolutely continuous by inspecting each term in \eqref{energy definition}.

From Proposition \ref{propexistence}, it follows that $\frac{1}{2}\norm{y(t)-\bar{y}}^{2}$ is absolutely continuous in $[0,+\infty)$.

In addition, for a.e. $t\geq 0$, it holds \begin{equation*}
    Ax(t)-f=Ax(t)-f^{\delta} + f^{\delta}-f=\dot{y}(t)+f^{\delta}-f 
\end{equation*}
 and since $\dot{y}$ is bounded in $[0,+\infty)$ (see point $2.$ in Proposition \ref{propexistence}), it follows that the function $t\norm{Ax(t)-f}^{2}$ is locally bounded and thus locally integrable.
 Therefore, the functions $\int_{0}^{t}s\norm{Ax(s)-f}^{2} ds$, $\int_{0}^{t}\norm{y(s)-\bar{y}} ds$ and  $\int_{0}^{t}s\norm{Ax(s)-f}  ds$ are absolutely continuous in $[0,+\infty)$.

Finally, since $-A^{\ast}\bar{y}\in\partial J(\bar{x})$and $-A^{\ast}y(t)\in\partial J(x(t))$ by the Fenchel-Moreau inequality,(see e.g. \cite[Theorem $16.29$]{bauschke2011convex}), it holds:
\begin{align}
J(x(t)) +J^{\ast}(-A^{\ast}y(t)) & = -\scal{A^{\ast}y(t)}{x(t)} \label{aligndual1} \\
\text{and } \quad J(\bar{x}) +J^{\ast}(-A^{\ast}\bar{y}) &= -\scal{A^{\ast}\bar{y}}{\bar{x}} \label{aligndual2}
\end{align}
By using \eqref{aligndual1} and \eqref{aligndual2} in the definition of $D^{-A^{\ast}y(t)}\left(\bar{x},x(t)\right)$, we find:
\begin{equation}\label{bregmanAC}
    \begin{aligned}
        D^{-A^{\ast}y(t)}(\bar{x}&,x(t))=J(\bar{x})-J(x(t)) +\scal{A^{\ast}y(t)}{\bar{x}-x(t)} \\
        & = J^{\ast}(-A^{\ast}y(t))-J^{\ast}(-A^{\ast}\bar{y}) +\scal{A^{\ast}y(t)}{x(t)} \\ & \quad  -\scal{A^{\ast}\bar{y}}{\bar{x}} +\scal{A^{\ast}y(t)}{\bar{x}-x(t)} \\ & =J^{\ast}(-A^{\ast}y(t))-J^{\ast}(-A^{\ast}\bar{y}) \\ & \quad + \scal{A^{\ast}y(t)-A^{\ast}\bar{y}}{\bar{x}} \\
     & =J^{\ast}(-A^{\ast}y(t))-J^{\ast}(-A^{\ast}\bar{y}) +\scal{f^{\delta}}{y-\bar{y}} \\
     & \quad +\scal{f-f^{\delta}}{y-\bar{y}} \\
     &= \mathcal{D}_{\delta}(y(t))-\mathcal{D}_{\delta}(\bar{y})+\scal{f-f^{\delta}}{y-\bar{y}}
    \end{aligned}
\end{equation}
From \eqref{bregmanAC} and point $4.$ of Proposition \ref{propexistence}, it follows that the function $D^{-A^{\ast}y(t)}\left(\bar{x},x(t)\right)$ is absolutely continuous with:
\begin{equation}\label{derivativebregman}
    \begin{aligned}
        \frac{d}{dt}D^{-A^{\ast}y(t)}\left(\bar{x},x(t)\right) &= \frac{d}{dt}\mathcal{D}_{\delta}(y(t)) +\scal{\dot{y}(t)}{f-f^{\delta}} \\
        &= -\norm{Ax(t)-f^{\delta}}^{2} +\scal{\dot{y}(t)}{f-f^{\delta}}
    \end{aligned}
\end{equation}

In a similar way, by using \eqref{aligndual1} and \eqref{aligndual2}, for a.e. $t\geq 0$, it holds: 
\begin{equation}\label{bregmanAC2}
    \begin{aligned}
        D^{-A^{\ast}\bar{y}}(x(t),\bar{x}) &= J^{\ast}(-A^{\ast}\bar{y}) +\scal{f^{\delta}}{\bar{y}} - J^{\ast}(-A^{\ast}y(t)) \\
        & \quad  -\scal{f^{\delta}}{y(t)} -\scal{Ax(t)-f^{\delta}}{y(t)-\bar{y}}\\
        & =\mathcal{D}_{\delta}(\bar{y}) -\mathcal{D}_{\delta}(y(t))-\scal{\dot{y}(t)}{y-\bar{y}} 
    \end{aligned}
\end{equation}
By Proposition \ref{propexistence}, the right-hand side of \eqref{bregmanAC2} is locally integrable, therefore it follows that for all $t\geq 0$, the function $\int_{0}^{t}D^{-A^{\ast}\bar{y}}(x(s),\bar{x}) ds$ is absolutely continuous.

Therefore, it follows that the function $V:[0,+\infty)\to\R$ is absolutely continuous as sum of absolutely continuous functions.

By taking the derivative of $V$ with respect to time, and using \eqref{derivative y} and \eqref{derivativebregman}, for a.e. $t\geq 0$ we obtain:
\begin{equation}\label{eq: derivativeV}
    \begin{aligned}
        \dot{V}(t)&= -\left(D^{-A^{\ast}\bar{y}}(x(t),\bar{x})  +D^{-A^{\ast}y(t)}(\bar{x},x(t))\right) \\ &+\langle f-f^{\delta}, y(t)-\bar{y}\rangle+t\norm{Ax(t)-f}^{2}\\
        & \quad +D^{-A^{\ast}\bar{y}}(x(t),\bar{x})   \\
        & \quad -\delta\norm{y(t)-\bar{y}} -\delta t\norm{Ax(t)-f} \\
        & \quad + D^{-A^{\ast}y(t)}\left(\bar{x},x(t)\right) \\
        & \quad -t\norm{Ax(t)-f^{\delta}}^{2} +t\scal{\dot{y}(t)}{f-f^{\delta}} \\
        &= \langle f-f^{\delta}, y(t)-\bar{y}\rangle  -\delta\norm{y(t)-\bar{y}} \\
        & \quad +
        t\norm{Ax(t)-f}^{2} -t\norm{Ax(t)-f^{\delta}}^{2} 
        \\
        & -\delta t\norm{Ax(t)-f} +t\scal{\dot{y}(t)}{f-f^{\delta}} 
    \end{aligned}
\end{equation}

By the KKT conditions \eqref{DGFP}, it holds:
\begin{equation}\label{eq: for derivativeV}
    \begin{aligned}
        \scal{\dot{y}(t)}{f-f^{\delta}} & = \scal{Ax(t)-f^{\delta}}{f-f^{\delta}} \\
         & =  -\norm{Ax(t)-f}^{2} +\norm{Ax(t)-f^{\delta}}^{2} \\
        &\quad + \scal{f-f^{\delta}}{f-A(x)}
    \end{aligned}
\end{equation}

By injecting relation \eqref{eq: for derivativeV} in \eqref{eq: derivativeV}, it follows that
\begin{equation}\label{eq:derivativeV2}
    \begin{aligned}
        \dot{V}(t) & = \langle f-f^{\delta}, y(t)-\bar{y}\rangle  -\delta\norm{y(t)-\bar{y}} \\
        & \quad + t\scal{f-f^{\delta}}{f-A(x)} -\delta t\norm{Ax(t)-f}
        \end{aligned}
\end{equation}
Finally by using Cauchy-Schwarz inequality for the scalar products in \eqref{eq:derivativeV2}, and the assumption on the noise \eqref{assumptiondelta}, it follows that
\begin{equation}
    \begin{aligned}
        \dot{V}(t) & \leq \norm{f-f^{\delta}}\norm{y(t)-\bar{y}}  -\delta\norm{y(t)-\bar{y}} \\
       ~ &  + t\norm{f-f^{\delta}}\norm{A(x)-f} -\delta t\norm{Ax(t)-f} \leq 0
    \end{aligned}
\end{equation}
which allows to conclude the proof of Proposition \ref{propVdecreasing}.

  \qed

\subsection{Proof of Theorem \ref{basiccor}}\label{subsectionProofteo}

We now provide the proof of the main Theorem \ref{basiccor}.

{\em Proof.} Since $D^{-A^{\ast}\bar{y}}(x(t),\bar{x})$ is non-negative and convex with respect to $x(t)$, by considering $\hat{x}(t)=\frac{2}{t}\int_{\frac{t}{2}}^{t}x(s) ds$ and using Jensen's inequality, for all $t\geq 0$, it holds: \begin{equation*}
    \begin{aligned}
        D^{-A^{\ast}\bar{y}}(\hat{x}(t),\bar{x}) & \leq \frac{2}{t}\int_{\frac{t}{2}}^{t}D^{-A^{\ast}\bar{y}}(x(s),\bar{x}) ds \\ & \leq \frac{2}{t}\int_{0}^{t}D^{-A^{\ast}\bar{y}}(x(s),\bar{x}) ds 
    \end{aligned}
    \end{equation*}
Therefore, the bound in \eqref{rateBregman} follows from the estimate \eqref{boundDsym} in Lemma \ref{lemma bound y}.

By using the non-increasing property of $V$ from Proposition \ref{propVdecreasing}, for all $t\geq 0$, we have:
\begin{equation}\label{eq: gros lyap beta}
    \begin{aligned}
       & \frac{1}{2}\norm{y(t)-\bar{y}}^{2}  +tD^{-A^{\ast}y(t)}\left(\bar{x},x(t)\right) \\ &+\int_{0}^{t} D^{-A^{\ast}\bar{y}}(x(s),\bar{x}) ds +\int_{0}^{t}s\norm{Ax(s)-f}^{2}ds \leq \\ &~~   \leq  V(0) + \delta\int_{0}^{t}\norm{y(s)-\bar{y}} ds +\delta\int_{0}^{t}s\norm{Ax(s)-f}  ds
    \end{aligned}
\end{equation}
By using the estimation \eqref{boundy} in Lemma \ref{lemma bound y} and setting $C_{0}=\norm{y(0)-\bar{y}}=2\sqrt{V(0)}$, we find:
\begin{equation}\label{eq: gros lyap}
    \begin{aligned}
       & \frac{1}{2}\norm{y(t)-\bar{y}}^{2}  +tD^{-A^{\ast}y(t)}\left(\bar{x},x(t)\right) \\ &+\int_{0}^{t} D^{-A^{\ast}\bar{y}}(x(s),\bar{x}) ds +\int_{0}^{t}s\norm{Ax(s)-f}^{2}ds \leq  \\ &~~   \leq  \frac{C_{0}^{2}}{2} +C_{0}\delta t +\frac{\delta^{2}t^{2}}{2}+\delta\int_{0}^{t}s\norm{Ax(s)-f}  ds
    \end{aligned}
\end{equation}
By using the Cauchy-Schwarz inequality for $\int_{0}^{t} s \norm{Ax(s)-f} ds$ in \eqref{eq: gros lyap}, it holds:
\begin{equation}\label{eq: gros lyap2}
    \begin{aligned}
       & \frac{1}{2}\norm{y(t)-\bar{y}}^{2}  +tD^{-A^{\ast}y(t)}\left(\bar{x},x(t)\right) \\ & \quad +\int_{0}^{t} D^{-A^{\ast}\bar{y}}(x(s),\bar{x}) ds +\int_{0}^{t}s\norm{Ax(s)-f}^{2}ds \leq  \\ &  \leq  \frac{C_{0}^{2}}{2} +C_{0}\delta t +\frac{\delta^{2}t^{2}}{2} \\
       &~~ \quad + \delta\left(\int_{0}^{t}s ds\right)^{\frac{1}{2}} \left(\int_{0}^{t} s \norm{Ax(s)-f}^{2} ds\right)^{\frac{1}{2}} \\
       & =\frac{(C_{0}+\delta t)^2}{2} +\frac{\delta t}{\sqrt{2}}\left(\int_{0}^{t} s \norm{Ax(s)-f}^{2} ds\right)^{\frac{1}{2}}
    \end{aligned}
\end{equation}
In particular, by neglecting the first three non-negative terms in the left-hand-side of \eqref{eq: gros lyap2}, we have: 
\begin{equation}\label{eq: sqrt brezis}
\begin{aligned}
    \int_{0}^{t}s\norm{Ax(t)-f}^{2}  ds & \leq \frac{(C_{0}+\delta t)^2}{2} \\ & \quad +\frac{\delta t}{\sqrt{2}}\left(\int_{0}^{t} s \norm{Ax(s)-f}^{2} ds\right)^{\frac{1}{2}}
\end{aligned}
\end{equation}

From \eqref{eq: sqrt brezis}, it follows that for every $t\geq 0$, it holds:
\begin{equation}\label{eq: totakesquares}
\begin{aligned}
   \left(\int_{0}^{t} s \norm{Ax(s)-f}^{2} ds\right)^{\frac{1}{2}} & \leq \frac{C_{0}+2\delta t}{\sqrt{2}}
\end{aligned}
\end{equation}

By plugging estimation \eqref{eq: totakesquares} into \eqref{eq: gros lyap2}, we find:
\begin{equation}\label{teleytaio}
\begin{aligned}
  tD^{-A^{\ast}y(t)}\left(\bar{x},x(t)\right) & +   \int_{0}^{t}s\norm{Ax(t)-f}^{2}  ds \leq 
 \\
 & \leq \frac{C_{0}^{2}}{2}+\frac{3C_{0}}{2}\delta t
 +\frac{3}{2}(\delta t)^{2}
\end{aligned}
 \end{equation}
which allows to derive the bound \eqref{estimateDinverse}.

Finally by using Jensen's inequality for $\hat{x}(t)=\frac{2}{t}\int_{\frac{t}{2}}^{t}x(s) ds$ in \eqref{teleytaio} and neglecting the non-negative term $tD^{-A^{\ast}y(t)}\left(\bar{x},x(t)\right)$, we obtain:
\begin{equation}
    \begin{aligned}
        & \norm{A\hat{x}(t)-f}^{2} \leq \frac{2}{t}\int_{\frac{t}{2}}^{t}\norm{Ax(s)-f}^{2} ds \\
        & 
\quad \leq \frac{4}{t^{2}}\int_{\frac{t}{2}}^{t}s \norm{Ax(s)-f}^{2} ds \leq \frac{4}{t^{2}}\int_{0}^{t}s\norm{Ax(t)-f}^{2}  ds \\ & \quad  \leq \frac{2C_{0}^{2}}{t^{2}}+6C_{0}\frac{\delta}{t} +6\delta^{2},
    \end{aligned}
    \end{equation}
    which allows to derive \eqref{ratefeasibility} and conclude the proof of Theorem \ref{basiccor}.
  \qed

\section{Conclusion and perspectives}\label{sectionconcl}
In this work we have studied the implicit regularization properties of the dual subgradient flow induced by a convex bias. As already commented in Paragraph \ref{secdiscr}, it would be interesting to analyze analogous properties for the algorithm proposed in \eqref{subgradalg}. Another interesting future direction is to extend the stability results to second order in-time dynamical systems and algorithms, which may provide the same stability results with earlier stopping times (and thus less computational resources) (see e.g. \cite{apidopoulos2018differential,boct2018second}).
Finally, a future line of research is related with the study of regularization properties of the continuous flow associated to the primal-dual algorithm \cite{chambolle2011first,molinari2021iterative}. 

\section*{APPENDIX}
The following lemma is a particular version of Bihari's lemma (see \cite{bihari1956generalization}). Its proof can be found in \cite[Lemma A.$5$]{brezis1973ope}. 
\begin{lemma}\label{lemma Bihari}
	Let $a$ be a non-negative constant; 
 and $\psi(t)$, $k(t)$ be positive and continous functions in $[0,T]$, such that the following inequality is satisfied:
	\begin{equation}\label{hyp}
	\psi(t)^{2} \leq a+\int_{t_{0}}^{t} k(s)\psi(s) \ ds \ \ \ \ \text{for every} \ t_{0}\leq t \leq t_{1}
	\end{equation}
Then for all $t\in [0,T]$, it holds:
	$$\psi(t) \leq \sqrt{a} +\frac{1}{2}\int_{0}^t k(s) ds$$
\end{lemma}

\subsection{Proofs for Tikhonov regularization path}\label{TYKgen}

Here we provide the proof of Theorem \ref{teotikh} concerning the Tikhonov regularization path: 
\begin{equation}\label{TRPP}
\begin{cases}
u \in & \partial J^\star (-sA^*v);\\
v=& Au-f^{\delta}.
\end{cases}
\end{equation}
By using the definition of $D^{-A^{\ast}\bar{y}}(u(s),\bar{x})$ and  $D^{-sA^{\ast}v(s)}(\bar{x},u(s))$ and the system \eqref{TRPP}, it holds:
\begin{equation}\label{jkfa}
\begin{split}
D^{-A^{\ast}\bar{y}}&(u(s),\bar{x})  +D^{-sA^{\ast}v(s)}(\bar{x},u(s))  = \\
& = -s \norm{Au(s)-f}^2  + \langle s\left(f^{\delta}-f\right)+\ys, Au(s)-f\rangle\\
& \leq -s \norm{Au(s)-f}^2 + \frac{1}{2\varepsilon}\norm{s\left(f^{\delta}-f\right)+\ys}^2 \\& \quad +\frac{\varepsilon}{2}\norm{Au(s)-f}^2\\
& \leq -\left(s-\frac{\varepsilon}{2}\right) \norm{Au(s)-f}^2 + \frac{1}{2\varepsilon}\left(s\delta+\norm{\ys}\right)^{2}
\end{split}
\end{equation}
where in the first inequality we used Young's inequality with $\varepsilon>0$ and in the second inequality the Cauchy-Schwarz inequality for the product $\scal{s\left(f^{\delta}-f\right)}{\bar{y}}$.

By setting $\varepsilon =2\eta s$, 
 for any $\eta\in(0,1]$ in \eqref{jkfa}, for all $s>0$, we find:
\begin{equation}
\begin{aligned}
    D^{-A^{\ast}\bar{y}}&(u(s),\bar{x})  +D^{-sA^{\ast}v(s)}(\bar{x},u(s)) ~+ \\
   & +(1-\eta)s\norm{Au(s)-f}^{2}  \leq \frac{1}{4\eta s}\left(\norm{\bar{y}}+s\delta\right)^2
    \end{aligned}
\end{equation}
Finally by setting $\eta=1$ and $\eta=\frac{1}{2}$ we deduce \eqref{DsymTikh} and \eqref{feasibilityTykh} respectively.


\section*{ACKNOWLEDGMENT}

We acknowledge the financial support of the European
Research Council (grant SLING 819789), the AFOSR ((European Office of Aerospace
Research and Development)) project FA9550-18-1-7009 and FA8655-22-1-7034, the EU H2020-MSCA-RISE project NoMADS - DLV-777826, the H2020-MSCA-ITN Project Trade-OPT 2019;  L. R. acknowledges  the Center
for Brains, Minds and Machines (CBMM), funded by NSF STC award CCF-1231216 and IIT. 
S.V. and C.M. are part of the Indam group ``Gruppo Nazionale per
l'Analisi Matematica, la Probabilit\`a e le loro applicazioni''.


\bibliographystyle{plain}
\bibliography{reference}

\end{document}